\documentclass[12pt]{article}
\begin{document}
\author{M.F.Hoque\\Department of Mathematics\\Pabna University of Science and Technology, Bangladesh\\E-mail: fazlul\_math@yahoo.co.in\\\\A.C.Paul\\Department of Mathematics\\Rajshahi University, Bangladesh\\Email: acpaulrubd\_math@yahoo.com} 

\title{Prime Gamma Rings with Centralizing and Commuting Generalized Derivations}  
\maketitle

\newtheorem{th1}{Theorem}[section]
\newtheorem{df}{Definition}
\newtheorem{lem}{Lemma}[section]
\newtheorem{cor}{Corollary}[section]
\newtheorem{pro}{Proposition}[section] 
\newtheorem{rem}{Remark}[section]

\begin{abstract} Let $M$ be a prime $\Gamma$-ring satisfying a certain assumption and $D$ a nonzero derivation on $M$. Let $f:M\rightarrow M$ be a generalized derivation such that $f$ is centralizing and commuting on a left ideal $J$ of $M$. Then we prove that $M$ is commutative.
\end{abstract}

\section{Introduction}An extensive generalized concept of classical ring set forth the notion of a gamma ring theory. As an emerging field of research, the research work of classical ring theory to 
the gamma ring theory has been drawn interest of many algebraists and prominent mathematicians over the world to determine many basic properties of gamma ring and to enrich the world of algebra.
The different researchers on this field have been doing a significant contributions to this field from its inception. In recent years, a large number of researchers are engaged to increase the efficacy of the results of gamma ring theory over the world.
\\\\The concept of a $\Gamma $-ring was first introduced by Nobusawa[11] and also shown that $\Gamma $-rings, more general than rings. Bernes[1] weakened slightly the conditions in the definition of $\Gamma $-ring in the sense of Nobusawa. 
They obtain a large numbers of improtant basic properties of $\Gamma $-rings in various ways and determined some more remarkable results of $\Gamma $-rings.
We start with the following necessary definitions.
\\\\ Let $M$ and $\Gamma $  be additive abelian groups. If there exists a mapping $(x,\alpha ,y)\rightarrow  x\alpha y$ of $M\times\Gamma \times M\rightarrow M$, which satisfies the conditions \\(i) $x\alpha y\in M$\\(ii) $(x+y)\alpha z$=$x\alpha z$+$y\alpha z$, $x(\alpha +\beta )z$=$x\alpha z$+$x\beta z$, $x\alpha (y+z)$=$x\alpha y$+$x\alpha z$\\(iii) $(x\alpha y)\beta z$=$x\alpha (y\beta z)$ for all $x,y,z\in M$ and $\alpha ,\beta \in \Gamma$, \\then $M$ is called a $\Gamma$-ring.
\\Every ring $M$ is a $\Gamma $-ring with $M$=$\Gamma$. However a $\Gamma $-ring need not be a ring.
\\Let $M$ be a $\Gamma $-ring. Then an additive subgroup $U$ of $M$ is called a left (right) ideal of $M$ if $M\Gamma U\subset U$($U\Gamma M \subset U$). If $U$ is both a left and a right ideal , then we say $U$ is an ideal of $M$. 
Suppose again that $M$ is a $\Gamma $-ring. Then $M$ is said to be a 2-torsion free if $2x$=$0$ implies $x$=$0$ for all $x\in M$.
An ideal $P_{1}$ of a $\Gamma $-ring $M$ is said to be prime if for any ideals $A$ and $B$ of $M$, $A\Gamma B\subseteq P_{1}$ implies $A\subseteq P_{1}$ or $B\subseteq P_{1}$. An ideal $P_{2}$ of a $\Gamma $-ring $M$ is said to be semiprime if for any ideal $U$ of $M$, $U\Gamma U\subseteq P_{2}$ implies $U\subseteq P_{2}$. A $\Gamma$-ring $M$ is said to be prime if $a\Gamma M\Gamma b$=$(0)$ with $a,b\in M$, implies $a$=$0$ or $b$=$0$ and semiprime if $a\Gamma M\Gamma a$=$(0)$ with $a\in M$ implies $a$=$0$. Furthermore, $M$ is said to be commutative $\Gamma $-ring if $x\alpha y$=$y\alpha x$ for all $x, y\in M$ and $\alpha \in\Gamma $. Moreover,the set $Z(M)$ =$\{x\in M:x\alpha y=y\alpha x $ for all $ \alpha \in \Gamma,  y\in M\}$ is called the centre of the $\Gamma $-ring $M$.
\\If $M$ is a $\Gamma $-ring, then $[x,y]_{\alpha }$=$x\alpha y-y\alpha x$ is known as the commutator of $x$ and $y$ with respect to $\alpha $, where $x,y\in M$ and $\alpha \in\Gamma $. We make the basic commutator identities:
\\  $[x\alpha y,z]_{\beta }$=$[x,z]_{\beta }\alpha y+x[\alpha ,\beta ]_{z}y+x\alpha [y,z]_{\beta }$ \\and  $[x,y\alpha z]_{\beta }$=$[x,y]_{\beta }\alpha z+y[\alpha ,\beta ]_{x}z+y\alpha [x,z]_{\beta }$ , \\for all $x,y.z\in M$ and $\alpha ,\beta \in\Gamma $.
We consider the following assumption:\\$(A)$.................$x\alpha y\beta z$=$x\beta y\alpha z$, for all $x,y,z\in M$, and $\alpha ,\beta \in\Gamma $.\\ According to the assumption $(A)$, the above two identites reduce to \\ $[x\alpha y,z]_{\beta }$=$[x,z]_{\beta }\alpha y+x\alpha [y,z]_{\beta }$ \\and  $[x,y\alpha z]_{\beta }$=$[x,y]_{\beta }\alpha z+y\alpha [x,z]_{\beta }$, \\which we extensively used.
\\\\ Note that Borut Zalar [14] worked on centralizers of semiprime rings and proved that Jordan centralizers and centralizers of this rings coincide. Joso Vukman[12, 13] developed some remarkable results using centralizers on  prime and semiprime rings.
Bresar[2], Mayne[10] and J.Luh[9] have developed some remarkable result on prime rings with commuting and centralizing.
Bernes[1], Luh [8] and Kyuno[7] studied the structure of $\Gamma$-rings and obtained various generalizations of corresponding parts in ring theory.
Y.Ceven [3] worked on Jordan left derivations on completely prime $\Gamma$-rings. He investigated the existence of a nonzero Jordan left derivation on a completely prime $\Gamma$-ring that makes the $\Gamma$-ring commutative with an assumption. With the same assumption, he showed that every Jordan left derivation on a completely prime $\Gamma$-ring is a left derivation on it.
\\\\An additive mapping $D:M\rightarrow M$ is called a derivation if $D(x\alpha y)=D(x)\alpha y+x\alpha D(y)$ holds for all $x,y\in M$, and $\alpha \in\Gamma $. 
A mapping $f$ is said to be commuting on a left ideal $J$ of $M$ if $[f(x),x]_{\alpha }=0$ for all $x\in J$, $\alpha \in\Gamma $ and $f$ is said to be centralizing if $[f(x),x]_{\alpha }\in Z(M)$ for all $x\in J$, $\alpha \in\Gamma $. 
An additive mapping $f:M\rightarrow M$ is said to be a generalized derivation on $M$ if $f(x\alpha y)=f(x)\alpha y+x\alpha D(y)$ for all $x,y\in M$ and $\alpha \in\Gamma $, $D$ a derivation on $M$.
\\\\In [3], M.F. Hoque and A.C Paul have proved that every Jordan centralizer of a 2-torsion free semiprime $\Gamma$-ring is a centralizer. There they also gave an example of a Jordan centralizer which is not a centralizer.
In [4], the same authors have proved that if $M$ is a 2-torsion free semiprime $\Gamma$-ring satisfying the assumption (A) and if $T:M\rightarrow M$ is an additive mapping such that $T(x\alpha y\beta x)=x\alpha T(y)\beta x$ for all $x, y\in M$ and $\alpha, \beta \in\Gamma$, then $T$ is a centralizer. Also, they have proved that $T$ is a centralizer if $M$ contains a multiplicative identity 1. 
Again in [5], they have proved that if $M$ be a 2-torsion free semiprime $\Gamma$-ring satisfying the assumption (A) and let  $T:M\rightarrow M$ be an additive mapping such that $2T(a\alpha b\beta a)=T(a)\alpha b\beta a+a\alpha b\beta T(a)$ holds for all pairs $a,b\in M$, and $\alpha , \beta \in\Gamma $. Then $T$ is a centralizer.
\\\\In this article, we prove that a prime $\Gamma $-ring $M$ is commutative if $f$ is a generalized derivation on $M$ with an assiciated nonzero derivation $D$ on $M$ such that $f$ is centralizing and commuting on a left ideal $J$ of $M$.

\vskip.2cm
\noindent
\section{Some priliminary results}
\vskip.2cm
\noindent  
We have to make some use of the following well-known results.
\begin{rem} Let $M$ be a prime $\Gamma $-ring. If $a\alpha b\in Z(M)$ with $0\neq a\in Z(M)$, then $b\in Z(M).$ 
\end{rem}
\vskip.2cm
\noindent 
\begin{rem} Let $M$ be a prime $\Gamma $-ring and $J$ a nonzero left ideal of $M$. If $D$ is a nonzero derivation on $M$, then $D$ is also a nonzero on $J$.
\end{rem}
\vskip.2cm
\noindent 
\begin{rem}Let $M$ be a prime $\Gamma $-ring and $J$ a nonzero left ideal of $M$. If $J$ is commutative, then $M$ is also commutative.
\end{rem}
We can easily prove the above remarks. Now we have to prove the following lemmas for the purpose of our main results.
\vskip.2cm
\noindent 
\begin{lem}Suppose $M$ is a prime $\Gamma$-ring satisfying the assumption (A) and $D:M\rightarrow M$ be a derivation. For an element $a\in M$, if $a\alpha D(x)=0$ for all $x\in M$ and $\alpha \in\Gamma $, then either $a=0$ or $D=0$.    
\end{lem}
\vskip.2cm
\noindent
{\bf Proof.} By our assumption, $a\alpha D(x)=0$ for all $x\in M$ and $\alpha \in \Gamma .$  We replace $x$ by $x\beta y$, then $a\alpha D(x\beta y)=0$ $\Rightarrow a\alpha D(x)\beta y+a\alpha y\beta D(x)=0$ $\Rightarrow a\alpha x\beta D(y)=0$ for all $x, y\in M$; $\alpha , \beta \in\Gamma $. 
If $D$ is not zero, that is, if $D(y)\neq 0$ for some $y\in M$, then by definition of prime $\Gamma $-ring, $a=0$. Hence proved. 
\vskip.2cm
\noindent
\begin{lem}Suppose $M$ is a prime $\Gamma$-ring satisfying the assumption (A) and $J$ a nonzero left ideal of M. If $M$ has a derivation $D$ which is zero on $J$, then $D$ is zero on $M$.
\end{lem}
\vskip.2cm
\noindent    
{\bf Proof.} By the hypothesis, $D(J)=0$. Replacing $J$ by $M\Gamma J$, we have $0=D(M\Gamma J)=D(M)\Gamma J+M\Gamma D(J)=D(M)\Gamma J$. Hence by Lemma-2.1, $D$ must be zero, since $J$ is nonzero.
\vskip.2cm
\noindent
\begin{lem}Suppose $M$ is a prime $\Gamma$-ring satisfying the assumption (A) and $J$ a nonzero left ideal of M. If $J$ is commutative on $M$, then $M$ is commutative.
\end{lem}
\vskip.2cm
\noindent  
{\bf Proof.} Suppose that $x$ is a fixed element in $J$. Since $J$ is commutative, so for all $y\in J$ and $\alpha \in\Gamma $, $[x,y]_{\alpha }=0$ and consequently, $[x,J]_{\Gamma }=0$.
Hence by Lemma-2.2, $[x,J]_{\Gamma }=0$ on $M$ and $x\in Z(M)$. Thus $[x,M]_{\Gamma }=0$ for every $x\in J$ and hence $[J,y]_{\Gamma }=0$ for all $y\in M$. Again Lemma-2.2, $[J,y]_{\Gamma }=0$ and $y\in Z(M)$ for all $y\in M$. Therefore $M$ is commutative.   
\vskip.2cm
\noindent   
\begin{lem}Let $M$ be a prime $\Gamma $-ring and $f:M\rightarrow M$ be an additive mapping. If $f$ is centralizing on a left ideal $J$ of $M$, then $f(a)\in Z(M)$ for all $a\in J\cup Z(M)$.
\end{lem}
\vskip.2cm   
\noindent
{\bf Proof.} By our assumption, $f$ is a centralizing on a left ideal $J$ of $M$. Thus we have, $[f(a),a]_{\alpha }\in Z(M)$ for all $a\in J$ and $\alpha \in\Gamma $. By linearization, for all $a, b\in J$ and $\alpha \in\Gamma $, we have \begin{eqnarray}&[f(a),b]_{\alpha }+[f(b),a]_{\alpha }\in Z(M)&\end{eqnarray} 
If $a\in Z(M)$, then equation (1) implies $[f(a),b]_{\alpha }\in Z(M)$. Now replacing $b$ by $f(a)\beta b$, we have $[f(a),f(a)\beta b]_{\alpha }\in Z(M)$, this implies $f(a)\beta [f(a),b]_{\alpha }\in Z(M).$
If $[f(a),b]_{\alpha }=0$, then $f(a)\in C_{\Gamma M}(J)$, the centralizer of $J$ in $M$ and hence $f(a)\in Z(M)$. Otherwise, if $[f(a),b]_{\alpha }\neq 0$, Remark-2.1 follows that $f(a)\in Z(M).$ Hence proved.
\vskip.2cm   
\noindent  
\section{The main results}
\vskip.2cm
\noindent  
\begin{th1}Let $M$ be a prime $\Gamma$-ring satisfying the assumption (A) and $D$ a nonzero derivation on $M$. If $f$ is a generalized derivation on a left ideal $J$ of $M$ such that $f$ is commuting on $J$, then $M$ is commutative.
\end{th1}         
\vskip.2cm 
\noindent  
{\bf Proof.} By our hypothesis, $f$ is commuting on $J$. Thus we have $[f(a), a]_{\alpha }=0$ for all $a\in J$ and $\alpha \in\Gamma$. By linearizing this relation implies $[f(a),b]_{\alpha }+[f(b),a]_{\alpha }=0.$
Putting $b=b\beta a$ and simplifying, we obtain $[b\beta D(a),a]_{\alpha }=0$. Replacing $b$ by $r\gamma b$, we have $[r,a]_{\alpha }\beta a\gamma D(a)=0$ for all $a\in J$, $r\in M$ and $\alpha ,\beta ,\gamma \in\Gamma $.
Since $M$ is prime $\Gamma $-ring, thus $[r,a]_{\alpha }=0$ or $D(a)=0.$ Therefore for any $a\in J$, either $a\in Z(M)$ or $D(a)=0$. Since $D$ is nonzero derivation on $M$,
then by Lemma-2.2, $D$ is nonzero on $J$. Suppose $D(a)\neq 0$ for some $a\in J$, then $a\in Z(M)$. Let $c\in J$ with $c\neq Z(M)$. Then $D(c)=0$  and $a+c\not\in Z(M)$, 
that is, $D(a+c)=0$ and so $D(a)=0$, which is a contradiction. Thus $c\in Z(M)$ for all $c\in J$. Hence $J$ is commutative and hence by Lemma-2.3, $M$ is commutative. 
\vskip.2cm
\noindent  
\begin{th1} Let $M$ be a prime $\Gamma $-ring satisfying the assumption (A) and $J$ a left ideal of $M$ with $J\cap Z(M)\neq 0$. If $f$ is a generalized derivation on $M$ with associated nonzero derivation $D$ such that $f$ is commuting on $J$, then $M$ is commutative.  
\end{th1}
\vskip.2cm 
\noindent 
{\bf Proof.} We claim that $Z(M)\neq 0$ because of $f$ is commuting on $J$ and the proof is complete. Now from equation (1), $[f(a),b]_{\alpha }+[f(b),a]_{\alpha }\in Z(M)$, if we replace $x$ by $c\beta b$ with $0\neq c\in Z(M)$,
then we have $[f(c),b]_{\alpha }\beta b+c\beta [D(b),b]_{\alpha }+c\beta [f(b),b]_{\alpha }\in Z(M)$. From Lemma-2.1, $f(c)\in Z(M)$ and hence $c\beta [D(b),b]_{\alpha }+c\beta [f(b),b]_{\alpha }\in Z(M)$. Since $f$ is a centralizing on $J$, we have $c\beta [f(b),b]_{\alpha }\in Z(M)$ and consequently $c\beta [D(b),b]_{\alpha }\in Z(M)$.  
As $c$ is nonzero, Remark-2.1 follows that $[D(b),b]_{\alpha }\in Z(M)$. This implies $D$ is centralizing on $J$ and hence we conclude that $M$ is commutative.   
 

\vskip.2cm 
\noindent             
\begin{thebibliography}{}
\bibitem{Eng1} W.E.Barnes, On the $\Gamma$-rings of Nobusawa, Pacific J.Math.,18(1966),411-422.
\bibitem{Eng2} M.Bresar, Centralizing Mappings and Derivations in Prime Rings, Journal of Algebra, 156(1993), 385-394.
\bibitem{Eng3} Y.Ceven, Jordan left derivations on completely prime gamma rings, C.U.Fen-Edebiyat Fakultesi,Fen Bilimleri Dergisi(2002)Cilt 23 Sayi 2.                                                                
\bibitem{Eng4} M.F.Hoque and A.C.Paul, On Centralizers of Semiprime Gamma Rings, International Mathematical Forum, Vol.6(2011), No.13,627-638.
\bibitem{Eng5} M.F.Hoque and A.C.Paul, Centralizers on Semiprime Gamma Rings, Italian J. Of Pure and Applied Mathematics, Vol. 30(2013), 289-302.
\bibitem{Eng6} M.F.Hoque and A.C.Paul, An Equation Related to Centralizers in Semiprime Gamma Rings, Annals of Pure and Applied Mathematics, Vol.1(1)(2012),84-90.
\bibitem{Eng7} S.Kyuno, On prime Gamma ring, Pacific J.Math.,75(1978),185-190. 
\bibitem{Eng8} J.Luh, A note on commuting automorphism of rings, Amer. Math. Monthly, 77(1970), 61-62.
\bibitem{Eng9} L.Luh, On the theory of simple Gamma rings, Michigan Math.J. 16(1969),65-75.
\bibitem{Eng10} J.Mayne, Centralizing automorphism of prime rings, Canad. Math. Bull. 19(1976), 113-115.
\bibitem{Eng11} N.Nobusawa, On the Generalization of the Ring Theory, Osaka J. Math.,1(1964),81-89.
\bibitem{Eng12} J.Vukman, Centralizers in prime and semiprime rings, Comment. Math. Univ. Carolinae 38(1997),231-240.
\bibitem{Eng13} J.Vukman, Centralizers on semiprime rings, Comment. Math. Univ. Carolinae 42,2(2001), 237-245. 
\bibitem{Eng14} B.Zalar, On centralizers of semiprime rings, Comment.Math. Univ. Carolinae 32(1991),609-614.  

\end{thebibliography}
\end{document}